\documentclass{amsart}
\usepackage{amssymb,latexsym,epsfig,color}

\newtheorem{theorem}{Theorem}[section]

\newtheorem{lemma}[theorem]{Lemma}

\theoremstyle{definition}

\def\proof{\smallskip\noindent {\it Proof: \ }}
\def\endproof{\hfill$\square$\medskip}
\def\R{\mathbb{R}}

\newcommand{\MON}{\mbox{\upshape \MON}}
\newcommand{\Int}{\mbox{\upshape int}\,}

\newcommand{\conv}{\mbox{\upshape conv}\,}
\newcommand{\rel}{\mbox{\upshape rel}\,}

\newcommand{\vol}{\mbox{\upshape vol}\,}
\newcommand{\vertex}{\mbox{\upshape vert}\,}

\title{How neighborly can a centrally symmetric polytope be?}

\author{Nathan Linial}
\author{Isabella Novik}
\address{Institute of Computer Science, Hebrew University, 
Jerusalem 91904, Israel, email: nati@cs.huji.ac.il}
\address{Department of Mathematics, University of
Washington, Box 354350, Seattle, WA 98195-4350, USA, email:
novik@math.washington.edu}    

\begin{document}

\begin{abstract}
We show that there exist
$k$-neighborly centrally symmetric $d$-dimensional polytopes
with $2(n+d)$ vertices, where
$$k(d,n)=\Theta\left(\frac{d}{1+\log ((d+n)/d)}\right).$$
We also show that this bound is tight.
\end{abstract}
\maketitle

\section{Introduction}
In this paper we study the following question:
{\bf How neighborly can a centrally symmetric polytope be
as a function of its dimension and the number of vertices?}
Let us recall the basic definitions: A polytope $P\subset \R^d$ is 
{\em centrally symmetric} (cs, for short) if for every $x\in P$, 
$-x$ belongs to $P$ as well. A cs polytope $P$ is called
{\em $k$-neighborly} if every set of $k$ of its vertices, no two
of which are antipodal is the vertex set of a face of $P$.

It is well-know that a general ({\sl non-cs}) $d$-dimensional
polytope  with at least $d+2$ vertices
can be at most $\lfloor d/2 \rfloor$-neighborly, and this bound 
is attained for instance by $d$-dimensional cyclic polytopes \cite[Example 0.6]{Ziegler}.
In contrast to the general case, the neighborliness of cs polytopes
appears to be quite restricted and not sufficiently understood.
A cs $d$-polytope with  at least $2(d+2)$ vertices
cannot be more than $\lfloor(d+1)/3\rfloor$-neighborly.
The case $d=4$ of this statement was observed by Gr\"unbaum \cite[p.116]{Gr}
in 1967, whereas the general case is due to McMullen and Shephard \cite{McMShep}.
These authors have also conjectured that a cs $d$-polytope with $2(d+n)$ vertices
cannot be more than 
   $\lfloor(d+n-1)/(n+1)\rfloor$-neighborly for all $n\geq 3$.
Their conjecture was refuted by Halsey \cite{Halsey}
and then by Schneider \cite{Schneider}, but only for $d>>n$. 
Namely, Schneider's theorem asserts that 
$$
\liminf_{d\rightarrow\infty} \frac{k(d,n)}{d+n} \geq 0.2390,
$$
where $k(d,n)$ denotes the largest
integer $k$ such that there exists a  $k$-neighborly cs $d$-polytope with $2(n+d)$ vertices.

A particularly interesting case is when $k(d,n)=1$.
Namely, given $d$, 
how large must $n$ be so that a cs $d$-polytope with $\geq 2(n+d)$ 
vertices cannot be even 2-neighborly?
That a cs $d$-polytope with a sufficiently large number of vertices 
($\approx (d/2)^{d/2}$) indeed cannot be 2-neighborly was verified by
Burton \cite{Burton}.
The McMullen-Shephard conjecture posits that this already holds for
$n = d-2$, but this turns out to be incorrect, since
we show that this critical $n$ is exponential in $d$.
No other facts on $k(d,n)$ appear to exist in the literature.

Here we compute the correct asymptotics of $k(d,n)$, thus
establishing an extension of 
Schneider's and  Burton's results.
Throughout the paper we set $m:=n+d$.

\begin{theorem} \label{main-thm}
$$\frac{C_1 d}{1+\log \frac{m}{d}} \leq
k(d,n) \leq 1+\frac{C_2 d}{1+\log \frac{m}{d}},
$$
where  $C_1, C_2>0$ are absolute constants independent of $d$ and $n$.
In particular, there exists a cs $d$-polytope 
with $4d$ vertices that is at least $\frac{d}{400}$-neighborly. 
\end{theorem}

Theorem \ref{main-thm} provides the correct asymptotic answer
for the above problem and shows that the largest
number of vertices in a 2-neighborly cs $d$-polytope
is $e^{\Theta(d)}$. In fact we can say a bit more here:

\begin{theorem} \label{main-thm2}
$k(d, 2^{d-1}+1-d)=1$. In other words,
a 2-neighborly cs $d$-polytope has at most
 $2^{d}$ vertices.
\end{theorem}

We prove Theorem \ref{main-thm2} as a warm-up for the proof of the upper bound
in Theorem~\ref{main-thm}. These proofs appear in Section 2.
The former result is a consequence of the Danzer-Gr\"unbaum theorem on 
the number of vertices of 
antipodal polytopes \cite{DanzGr} (a  more readily accessible source for a proof of this theorem
is Aigner and Ziegler's beautiful book \cite[Ch.~14]{AigZieg}) 
and the observation that every  2-neighborly cs
polytope is antipodal. For the latter result we establish a certain modification
of the Danzer-Gr\"unbaum argument.

The proof of the lower bound in 
Theorem~\ref{main-thm} is based on studying the 
cs transforms of cs polytopes introduced 
in \cite{McMShep} and on a theorem 
due to Garnaev and Gluskin \cite{GarnGlus}. This theorem concerns
the intersection of the $m$-dimensional octahedron $B_1^m$
with any $n$-dimensional subspace of $\R^m$. The question
is how close to an $n$-dimensional Euclidean ball
such an intersection can be.
We outline the necessary background on cs transforms in Section 3. 
The verification of the lower bound in Theorem \ref{main-thm} 
and the statement of the Garnaev-Gluskin theorem
are provided in Section~4.  The proof of Garnaev-Gluskin result
and hence  also of Theorem \ref{main-thm}  is probabilistic in nature: 
it does not give an explicit construction
of neighborly cs polytopes, but rather shows that they form a
set of positive probability in a certain probability space.
Indeed it is an interesting open question to find explicit constructions
of highly neighborly cs polytopes that meet the lower bound we prove.

\section{The upper bound on $k(d,n)$} The goal of this section is to verify 
Theorem \ref{main-thm2} and the upper bound in Theorem~\ref{main-thm}.
To this end, recall from \cite[p.~420]{Gr} that two vertices  
$v$ and $w$ of a  $d$-polytope $Q\subset\R^d$ 
are an {\em antipodal pair} if 
 there exist two parallel
(distinct) supporting hyperplanes of $Q$ one of which
contains $v$ and the other one contains $w$. A $d$-polytope $Q$
is called an {\em antipodal polytope} if every two of its vertices
form an antipodal pair.
The theorem due to Danzer and Gr\"unbaum 
\cite{DanzGr} (see also \cite[Ch.~14]{AigZieg}) 
asserts that an antipodal $d$-polytope
cannot have  more than $2^d$ vertices.  Theorem \ref{main-thm2}
is then an immediate consequence of their result and the following observation.

\begin{lemma}
Every  2-neighborly cs $d$-polytope is an antipodal
polytope.
\end{lemma}
\proof 
Consider two vertices $v$ and $w$ of a  2-neighborly cs $d$-polytope $P$. 
 We  show that
they form an antipodal pair. There are two possible cases:
either $v=-w$ or $v\neq -w$. In the first case, let $H$ be any 
hyperplane satisfying $H\cap P=\{v\}$. (Such an $H$ exists since $v$ is a vertex
of $P$). Then $H$ is a supporting hyperplane of $P$ that contains $v$,
while $-H:=\{x\in \R^d : -x\in H\}$ is a supporting hyperplane of $P$
that contains $w=-v$. In the second case, consider the set
$E=\conv\{v, -w\}$. 
Since $P$ is 2-neighborly,
$E$ is an edge of $P$, and so there exists
 a hyperplane $H$ such that $H\cap P=E$. Then $H$ 
is a supporting hyperplane of $P$ that contains $v$,
while $-H$ is a supporting hyperplane of $P$
that contains the edge $E'=-E=\conv \{-v, w\}$, and hence also 
the vertex $w$.
\endproof

The proof of the upper bound in 
Theorem \ref{main-thm} is obtained by a certain modification
of the Danzer-Gr\"unbaum argument and is
based on the following simple combinatorial result.
(We denote by ${[m] \choose s}$ the family of all $s$-subsets of $[m]:=\{1,\ldots, m\}$.)

\begin{lemma} \label{combin-lemma}
For every two
integers $s$ and $m$ with $s\leq m/2$, there exists a family 
$\mathcal{F}\subset {[m] \choose s}$ of size 
$\geq(Cm/s)^{ s/2}$ such that $|A\cap B|\leq s/2$
for every two distinct $A, B \in \mathcal{F}$.
Here $C > 0$ is an absolute constant.
\end{lemma}

\proof
We construct $\mathcal{F}$ by a simple greedy argument.
Let $\mathcal{F}\subset {[m] \choose s}$ be an inclusion-maximal family 
satisfying the condition $|A\cap B|\leq s/2$
for all $A, B\in \mathcal{F}$. We claim that
$|\mathcal{F}|\geq (Cm/s)^{ s/2 }$. Indeed,
for a set $A\in\mathcal{F}$ 
define the ``forbidden collection"
$A^f:=\{B\in {[m] \choose s} \ : \ |A\cap B|> s/2\}$.
Then 
$$|A^f|=\sum_{k=1}^{\lceil s/2 \rceil}
{s \choose \lfloor s/2\rfloor +k}
{m-s  \choose \lceil s/2\rceil-k} 
   < {m  \choose \lfloor s/2 \rfloor} 
\sum_{k=1}^{\lceil s/2 \rceil}
{s \choose \lfloor s/2\rfloor+k}
< 
 {m \choose \lfloor s/2 \rfloor}2^s,
$$ 
 and so 
$$|\bigcup_{A\in \mathcal{F}} A^f|\leq 
|\mathcal{F}| \cdot {m \choose \lfloor s/2 \rfloor}2^s \leq
|\mathcal{F}| \cdot 2^s\left(\frac{em}{s/3}\right)^{s/2} \leq
|\mathcal{F}| \cdot \left(\frac{12em}{s}\right)^{s/2}.$$
On the other hand, the maximality of $\mathcal{F}$
implies that 
$\mathcal{F} \cup (\bigcup_{A\in \mathcal{F}} A^f) =
{[m] \choose s}$. Hence
$|\mathcal{F}| \cdot(1+ (\frac{12em}{s})^
{s/2}) \geq
{m \choose s}$ which together with
the inequality ${m \choose s} \geq (\frac{m}{s})^s$
yields the result. \endproof

We are now in a position to prove the upper bound
in Theorem \ref{main-thm} asserting that
$k(d,n) \leq 1+\frac{C_2 d}{1+\log \frac{m}{d}}$
for some absolute constant $C_2>0$. To do so consider a cs
$d$-polytope $P$ on the vertex set
$V=\{v_1, -v_1, \ldots, v_m, -v_m\}$ that is $2s$-neighborly.
Thus $2s \leq d \leq m$.
Let $\mathcal{F} \subset {[m] \choose s}$ be a family from Lemma \ref{combin-lemma}.
For each set $A\in \mathcal{F}$, define
$$P_A:=P+\frac{2}{s}\sum_{i\in A} v_i \subseteq P+2P=3P$$
to be a translate of $P$, where ``+" denotes Minkowski addition. 

We claim that the polytopes $P_A$, $A\in\mathcal{F}$,
have pairwise disjoint interiors, whence
$$
|\mathcal{F}| \cdot \vol(P)=
\sum_{A\in \mathcal{F}} \vol(P_A)\leq \vol(3P)=3^d \cdot
\vol(P).$$ 
It follows that
$$
3^d \geq |\mathcal{F}|\geq \left(Cm/s\right)^
{s/2}.
$$
Thus 
$d\cdot \log 3 \geq  s/2  \cdot \log(Cm/s)
\geq \Omega(s \cdot \log(m/d))$, as claimed.

We turn to show that for any two distinct $A, B \in \mathcal{F}$
the sets $P_A$ and $P_B$ have disjoint interiors.
The symmetric difference $R:=(A-B)\cup (B-A)\subseteq [m]$ has
cardinality $s\leq |R| \leq 2s$
since $|A|=|B|=s$ and $|A\cap B|\leq s/2$.
Therefore, the fact that $P$ is $2s$-neighborly implies that
$\{-v_i: i\in A-B\} \cup \{v_j : j\in B-A\} \subset V$ is the vertex set,
$\vertex(F)$, of some proper face $F$ of $P$. We want to rule out
the possibility that
$x+\frac{2}{s} \sum_{i\in A}v_i = y+\frac{2}{s} \sum_{i\in B}v_i$
for some $x, y \in \Int(P)$. Indeed,
\begin{eqnarray*}
&& \frac{1}{2}\left[x+\frac{2}{s}\sum_{i\in A}v_i\right]-
 \frac{1}{2}\left[y+\frac{2}{s}\sum_{j\in B}v_j\right]\\&=&
\frac{x-y}{2} - \frac{1}{s}\left[\sum_{i\in A-B} (-v_i) +
 \sum_{j\in B-A} v_j\right]
=\frac{x-y}{2} - \frac{|R|}{s}\sum_{v\in 
\text{\tiny vert}(F)} \frac{1}{|R|} v\neq 0,
\end{eqnarray*}
since $(x-y)/2$ is an interior point of $P$,
while $\frac{|R|}{s}\sum_{v} \frac{1}{|R|} v$ 
is the $|R|/s \geq 1$-multiple of a boundary point of $P$, namely of the barycenter of $F$.
The conclusion follows.
\endproof

\section{Centrally Symmetric transforms}

Throughout the rest of the paper we denote the standard scalar product
on $\R^n$ by $\langle -,-\rangle$. The abbreviations int, rel int, and conv stand for
the interior, relative interior, and convex hull respectively. 

Following \cite{McMShep}, we define a
{\em centrally symmetric set} (cs set, for short)
as a finite spanning subset of $\R^d$
of the form $V=\{v_1, -v_1,\ldots,  v_m, -v_m\}$.
The construction described in 
\cite{McMShep} associates with a cs
set $V=\{\pm v_1, \ldots, \pm v_m\}\subset \R^d$ another
cs set $\overline{V}=
\{ \pm\overline{v}_1, \ldots, \pm\overline{v}_m \}\subset \R^{m-d}=\R^n$
called the {\em cs transform} of $V$. This operation
possesses the following properties.
\begin{enumerate}
\item Let $\overline{V}=
\{ \pm\overline{v}_1, \ldots, \pm\overline{v}_m \}\subset \R^{n}$ 
be a cs set. Then $\overline{V}$ is a cs transform of the vertex set $V$
 of a cs $d$-polytope $P$ with $2m$ vertices if and
only if 
$$
\overline{v}_i \in \Int \conv \left\{\sum_{l\in [m] -\{i\}}
\epsilon_l\overline{v}_l \; : \epsilon_l\in\{1, -1\}\right\} \quad 
\mbox{ for all  $i=1, \ldots, m$}.
$$

\item More generally, if $V=\{\pm v_1, \ldots, \pm v_m\}\subset \R^d$
 is the vertex set of a cs $d$-polytope $P$ and $\overline{V}$
is a cs transform of $V$, then the set
$\{\delta_i v_{i}: i\in I\}\subset V$ 
(where
 $\delta_i\in\{1,-1\}$, $i\in I$, are fixed signs 
and $I=\{i_1<\ldots<i_k\} \subseteq[m]$)
 is the vertex set of a face of $P$ if and only if
$$
\sum_{i\in I}\delta_i\overline{v}_{i}\in 
\rel \Int \conv \left\{\sum_{l\in [m] -I}
\epsilon_l \overline{v}_l \; : \epsilon_l\in\{1, -1\}\right\}.
$$
\end{enumerate}

Thus if $V=\{\pm v_1, \ldots, \pm v_m\}\subset \R^d$
 is a cs set, then its
subset  $\{\delta_i v_{i} : i \in I\}$ 
(for some $\delta_i\in\{1,-1\}$ and $I=\{i_1<\ldots<i_k\}\subseteq m$)
fails to be the vertex set of a face of the cs polytope
 $P:=\conv V$ if and only if
there exists $u\in \R^n$ 
such that
$$
0\neq \langle\sum_{i \in I}\delta_i\overline{v}_{i}, u\rangle \geq 
\langle\sum_{l\in [m] -I}
\pm \overline{v}_l, u\rangle \quad \mbox{ for all choices of signs},
$$
that is, if and only if
$$
0\neq \sum_{i \in I}\langle\delta_i\overline{v}_{i}, u\rangle \geq 
\sum_{l\in [m] -I}
|\langle \overline{v}_l, u\rangle|. 
$$

We call a subset $\{\overline{v}_{i} : i\in I\}$ of 
$\overline{V}_+:=\{\overline{v}_1, \ldots, \overline{v}_m\}\in\R^{n}$
{\em dominant} if there exists  $0\neq u\in\R^{n}$ such that
$$
\sum_{i\in I} |\langle \overline{v}_{i_j}, u\rangle|
\geq \sum_{l\in [m] -I} |\langle \overline{v}_l, u\rangle|, 
\mbox{ or equivalently, } 
\sum_{i \in I} |\langle \overline{v}_{i}, u\rangle|
\geq \frac{1}{2}\sum_{l=1}^m |\langle \overline{v}_l, u\rangle|.
$$

Since $\sum_{i\in I}\langle\delta_i\overline{v}_{i}, u\rangle \leq
\sum_{i\in I} |\langle \overline{v}_{i}, u\rangle|$ for $\delta_i=\pm 1$, and
since equality is attained  for a certain choice of signs, we obtain
the following criterion.
\begin{lemma}   \label{cs}
A cs set $\overline{V}
=\{\pm \overline{v}_1, \ldots, \pm \overline{v}_{m}\}\subset \R^{n}$ 
is a cs transform of the vertex set of a $k$-neighborly cs
$d$-polytope with $2m=2(n+d)$ vertices if and only if the set
$\overline{V}_+:=
\{\overline{v}_1, \ldots,  \overline{v}_{m}\}$ does not contain 
dominant subsets of size $k$.
\end{lemma}

Thus to prove lower bounds on $k(d,n)$ 
it suffices to construct  vector configurations  spanning $\R^n$ that do
not contain small dominant subsets. This is done in the following section.

\section{Vector configurations without small dominant subsets}
For a vector $x=(x_1, \ldots, x_m)\in\R^m$ its 
$l^m_1$ and $l^m_2$
norms are defined as
$\|x\|_1:=\sum_{i=1}^m |x_i|$  and 
$\|x\|_2:=\sqrt{\sum_{i=1}^m |x_i|^2}$,
respectively. Thus the unit ball of $l^m_1$ is the $m$-dimensional
octahedron $B^m_1:=\{x\in \R^m : \sum_{i=1}^m |x_i|\leq 1\}$,
while the unit ball of $l^m_2$
is the $m$-dimensional Euclidean ball 
$B^m_2=\{x\in \R^m : \sum_{i=1}^m x_i^2\leq 1\}$.
The theorem due to Garnaev and Gluskin \cite{GarnGlus} 
(see \cite{Makovoz}
for a simplified proof) quantifies the extent to which the intersection of 
$B^m_1$ with an $n$-dimensional  subspace of $\R^m$
can be close to $B^{n}_2$. It asserts that for any
natural numbers $d$ and $n$  
there exists a subspace $L^d$ of $\R^m$ of codimension $d$
(equivalently, of dimension $m-d=n$),
 such that 
\begin{equation}  \label{GG}
\|x\|_2 \leq \widetilde{C}
 \sqrt{\frac{1+\log(m/d)}{d}} \cdot\|x\|_1 \quad \mbox{ for all }
x\in L^d,
\end{equation}
where $\widetilde{C}$ is an absolute constant independent of $d$ and $m$. 
(In fact, the set of such subspaces has a positive measure
in the Grassmannian manifold $G_{m-d,m}$ of all codimension 
$d$ subspaces of $\R^m$ endowed with the normed unitary invariant measure.)
In the following we refer to such subspace $L^d$ as a Garnaev-Gluskin
subspace.
A weaker version of this theorem, with $(1+\log(m/d))^{3/2}$ 
instead of
$(1+\log(m/d))^{1/2}$, had been  shown
earlier by Ka\v sin \cite{Kashin}.

Since $L^d\subset \R^m$ is an $n$-dimensional space, 
there is a linear injective map  $T: \R^n \rightarrow \R^m$ 
whose image is $L^d$.
Let $A$ be the $m\times n$ matrix representing this map, and let
 $\overline{v}_1, \cdots, \overline{v}_m\in \R^n$ be the rows of this matrix.
Then for every $0\neq u\in\R^n$, $T(u)$ is a non-zero element of $L^d$
whose $i$-th  coordinate is  given by 
$\langle \overline{v}_i, u\rangle$ (for $i=1, \ldots, m$). 
Hence for every $k$-element subset
 $\{\overline{v}_{i} : i\in I\}$
of $\{\overline{v}_1, \cdots, \overline{v}_m\}$ and for every 
$0\neq u\in\R^n$, 
we have
\begin{eqnarray*}
\sum_{i\in I} |\langle \overline{v}_{i}, u\rangle| 
&\leq& \sqrt{k} \cdot
\sqrt{\sum_{i \in I} \langle \overline{v}_{i}, u\rangle^2}
\qquad \qquad \mbox{ (by the Cauchy-Schwarz inequality) }\\
&\leq& \sqrt{k} \cdot \|T(u)\|_2 \leq  \sqrt{k}\cdot \widetilde{C}\cdot
 \sqrt{\frac{1+\log(m/d)}{d}} \cdot\|T(u)\|_1
\quad \mbox{ (by Eq.~(\ref{GG}))} \\
&=&\sqrt{k}\cdot \widetilde{C}\cdot
\sqrt{\frac{1+\log(m/d)}{d}}\cdot \sum_{l=1}^m 
|\langle \overline{v}_{l}, u\rangle| 
< \frac{1}{2} \sum_{l=1}^m 
|\langle \overline{v}_{l}, u\rangle|
\end{eqnarray*}
as long as $k< \frac{1}{4\widetilde{C}^2}\cdot \frac{d}{1+\log(m/d)}$. 
Therefore
we infer the following result.
\begin{lemma}   \label{dominant}
Let $L^d$ be a Garnaev-Gluskin subspace of $\R^{m}$. Denote by 
 $\overline{V}_+=
\{\overline{v}_1, \cdots, \overline{v}_{m}\}\subset \R^n$ the 
set of rows of the matrix
representing a map $T: \R^n \rightarrow \R^{m}$ whose image is
$L^d$. Then $\overline{V}_+$ does not contain dominant
subsets of the size smaller than 
$\lceil \frac{1}{4\widetilde{C}^2}\cdot \frac{d}{1+\log\frac{m}{d}}\rceil$.
\end{lemma}

The lower bound in Theorem \ref{main-thm} asserting that for every $n$ and $d$
there exists a cs $d$-polytope with $2m=2(n+d)$ vertices that is
$\Omega(\frac{d}{1+\log \frac{m}{d}})$-neighborly
 is then an immediate corollary of Lemmas \ref{cs} and  \ref{dominant}.

To obtain an estimate on the constant for the $n=d$ case (the 
``in particular"-part of Theorem \ref{main-thm}), 
we use the following result essentially due to Ka\v sin
(see \cite[page 21]{Ball}) asserting that 
 there exists an orthogonal transformation  $U$ of $\R^d$ such that 
\begin{equation} \label{Szarek}
\|x\|_2 \leq \frac{4R^2}{\sqrt{d}}\left(\|U^{-1}x\|_1+\|x\|_1\right) 
\qquad \mbox{ for all }x\in\R^d.
\end{equation}
Here $R$ is the {\em volume ratio} of the octahedron $B^d_1$
(the notion introduced by Szarek \cite{Szarek}), that is,
\begin{equation} \label{vol-ratio}
R:=\left(\frac{\vol(\sqrt{d} B^d_1)}{\vol(B^d_2)}\right)^{1/d}
=\left( \frac{2^d d^{d/2}}{d!} \frac{\Gamma(d/2+1)}{\pi^{d/2}}\right)^
{1/d} \leq \left( \frac{2e}{\pi}\right)^{1/2},
\end{equation}
where  $\Gamma(\cdot)$ denotes the Gamma function.

Consider the $2d$-element set
$\overline{V}_+:=\{e_1, \ldots, e_d, Ue_1, \ldots, Ue_d\}\subset \R^d$.
The calculations completely analogous to those 
in the proof of Lemma \ref{dominant},
 but using Eq.~(\ref{Szarek})  instead of (\ref{GG}), imply
that all dominant subsets of the set $\overline{V}_+$ have size of at least
$$\frac{d}{2^7 R^4} \stackrel{\mbox{ {\small by }} (\ref{vol-ratio})}{\geq}
 \frac{\pi^2}{2^9 e^2} d \geq \frac{d}{400}. 
$$
This fact together with Lemma \ref{cs} yields the second part of the 
Theorem \ref{main-thm}, that is, existence of a
$\frac{d}{400}$-neighborly cs $d$-polytope with $4d$ vertices.

\section{Concluding remarks}
We close the paper with the following remark concerning
the exact value of $k(n,d)$.
Let $\overline{V}_+=\{\overline{v}_1, \ldots, \overline{v}_m\}\subset \R^n$
be an $m$-element set that spans $\R^n$.
 Denote by
 $A=A(\overline{V}_+)$ the $m\times n$ matrix
whose rows are the elements of $\overline{V}_+$. Let 
$T: \R^n \rightarrow \R^m$ be the map represented by $A$, and 
let $L^d=L^d(\overline{V}_+)$ be the image of $T$.
Clearly $L^d$ has codimension $d$ in $\R^m$, and 
every codimension $d$ subspace of $\R^m$ arises this way.

Fix an integer $s\in [m]$ and consider the
norm $|||-|||_s$ on $\R^m$ defined by 
$$|||(x_1, \ldots, x_m)|||_s:=\max\{ \sum_{i\in \sigma} |x_i| :
\sigma \subset [m], |\sigma|=s \}.$$
(E.g.~$|||-|||_1$ coincides with the  $l^m_\infty$ norm.) 
Since the $i$-th coordinate of $x=T(u)\in L_d$
equals $\langle \overline{v}_i, u \rangle$,
it follows that $\overline{V}_+$  does not have a dominant 
subset of size $s$ if and only if  
$|||x|||_s < \frac{1}{2}\|x\|_1$ for every $0\neq x\in L^d$. Thus Lemma~\ref{cs} implies that
\begin{equation} \label{Gelfand}
k(d,n)= \max \{s : c_d(l_1, |||-|||_s) <1/2 \} =
\min\{s: c_d(l_1, |||-|||_s) \geq 1/2 \} -1,
\end{equation}
where
$c_d=c_d(l_1, |||-|||_s)$ is the {\em $d$-th Gelfand number},
$$c_d:=\inf_{L^d\subset \R^m} \sup_{x\in L_d-\{0\}} \frac{|||x|||_s}{\|x\|_1}.  
$$
(Here the infimum is taken over all codimension $d$ subspaces of $\R^m$).
Gelfand numbers have received a good deal of attention in Banach Space
Theory \cite{Carl, Pinkus}.

It is a notoriously difficult question to find explicit constructions
for spaces that satisfy the conditions as in the work of Ka\v sin, 
Garnaev and Gluskin. Perhaps it is less difficult, though,
to construct highly neighborly cs polytopes. Any progress on this
problem would be of interest.


\begin{thebibliography}{999}
\bibitem{AigZieg} M.~Aigner and G.~M.~Ziegler,
{\em Proofs from THE BOOK}, Third edition, Springer-Verlag, Berlin, 2004.


\bibitem{Ball} K.~Ball,  {\em An elementary introduction to 
modern convex geometry}, 
Flavors of geometry, 1--58, Math. Sci. Res. Inst. Publ., {\bf 31}, 
Cambridge Univ. Press, Cambridge, 1997. 

\bibitem{Burton} G.~R.~Burton,  
``The nonneighbourliness of centrally symmetric convex polytopes 
having many vertices'', J.~Combin.~Theory Ser.~A {\bf 58} (1991), 321--322.

\bibitem{Carl}  B.~Carl and I.~Stephani,
{\em Entropy, compactness and the approximation of operators}, Cambridge Tracts in Mathematics, {\bf 98},
 Cambridge University Press, Cambridge, 1990.


\bibitem{DanzGr} L.~Danzer and B.~Gr\"unbaum, 
``\"Uber zwei Probleme bez\"uglich konvexer K\"orper von P.~Erd\"os 
und von V.~L.~Klee'' 
(German), Math.~Z. {\bf 79} (1962) 95--99.


\bibitem{GarnGlus}
 A.~Yu.~Garnaev and E.~D.~Gluskin, 
``The widths of a Euclidean ball'' (Russian),
 Dokl. Akad. Nauk SSSR {\bf 277} (1984), no. 5, 1048--1052.
(English translation: Soviet Math. Dokl. 30 (1984), no. 1, 200--204.)



\bibitem{Gr}
B.~Gr\"unbaum, {\em Convex polytopes}, Second edition 
(Prepared and with a preface by V.~Kaibel, V.~Klee and G.~M.~Ziegler),
 Graduate Texts in Mathematics, {\bf 221}, Springer-Verlag, New York, 2003.

\bibitem{Halsey} E.~R.~Halsey, ``Zonotopal complexes on the $d$-cube'',
{\em Doctoral dissertation, University of Washington} (1972).

\bibitem{Kashin} B.~S.~Ka\v sin, 
``The widths of certain finite-dimensional sets and 
classes of smooth functions'' (Russian),
 Izv.~Akad.~Nauk SSSR Ser.~Mat. {\bf 41} (1977), no. 2, 334--351, 478.
(English translation: Math. USSR-Izv. 11 (1977), no. 2, 317--333 (1978).)

\bibitem{Makovoz} 
Y.~Makovoz, 
``A simple proof of an inequality in the theory of $n$-widths'',
 Constructive theory of functions (Varna, 1987), 305--308, 
Publ. House Bulgar. Acad. Sci., Sofia, 1988. 

\bibitem{McMShep}
P.~McMullen, P. and G.~C.~ Shephard, 
 ``Diagrams for centrally symmetric polytopes'', Mathematika {\bf 15}
 (1968), 123--138.  

\bibitem{Pinkus}  A.~Pinkus, {\em On $L{\sp 1}$-approximation}, Cambridge Tracts in 
Mathematics, {\bf 93}, Cambridge University Press, Cambridge, 1989.


\bibitem{Schneider}
R.~Schneider, 
``Neighbourliness of centrally symmetric polytopes in high dimensions'',
Mathematika {\bf 22} (1975), no. 2, 176--181.

\bibitem{Szarek} S.~J.~Szarek, 
``On Kashin's almost Euclidean orthogonal decomposition of $l\sp{1}\sb{n}$", 
Bull. Acad. Polon. Sci. S\'er. Sci. Math. Astronom. Phys. {\bf 26} (1978),  
691--694.


\bibitem{Ziegler}  
G.~M.~Ziegler,  
{\em Lectures on polytopes},
 Graduate Texts in Mathematics, {\bf 152}, Springer-Verlag, New York, 1995. 

\end{thebibliography}
\end{document}